\documentclass[a4paper,12pt]{article} 
\usepackage{color,amssymb,amsthm,amsmath,amsxtra,multicol,nextpage,fullpage}

\usepackage{fullpage}

\newtheorem{theorem}{Theorem}

\newtheorem{corollary}{Corollary}

\newtheorem{lemma}{Lemma}

\newtheorem{proposition}[theorem]{Proposition}

\begin{document}

\title{Pieri-Type Formulas for the Nonsymmetric Macdonald Polynomials}
\author{W. Baratta \\
	Department of Mathematics, University of Melbourne}
	\maketitle

\begin{abstract}
In symmetric Macdonald polynomial theory the Pieri formula gives the
branching coefficients for the product of the $r^{th}$ elementary symmetric
function $e_{r}(z)$ and the Macdonald polynomial $P_{\kappa }\left( z\right) 
$. In this paper we give the nonsymmetric analogues for the cases $r=1$ and $%
r=n-1.$ We do this by first deducing the the decomposition for the product
of any nonsymmetric Macdonald polynomial $E_{\eta }\left( z\right) $ with $%
z_{i}$ in terms of nonsymmetric Macdonald polynomials. As a corollary of
finding the branching coefficients of $e_{1}\left( z\right) E_{\eta }\left(
z\right) $ we evaluate the generalised binomial coefficients $\tbinom{\eta }{%
\nu }$ associated with the nonsymmetric Macdonald polynomials for $%
\left\vert \eta \right\vert =\left\vert \nu \right\vert +1.$
\end{abstract}

\section{Introduction\label{introduction}}

The nonsymmetric Macdonald polynomials $E_{\eta }:=E_{\eta }\left(
z;q,t\right) $ are polynomials of $n$ variables $z=\left(
z_{1},...,z_{n}\right) $ having coefficients in the field $\mathbb{Q}
\left( q,t\right) $ of rational functions of the indeterminants $q$ and $t$.
The compositions $\eta :=\left( \eta _{1},...,\eta _{n}\right) $ of
non-negative integers parts $\eta _{i}$ label these polynomials. The
nonsymmetric Macdonald polynomials can be defined, up to normalisation, as
the unique simultaneous eigenfunctions of the commuting operators 
\begin{equation}
Y_{i}=t^{-n+1}T_{i}...T_{n-1}\omega T_{1}^{-1}...T_{i-1}^{-i},\text{ }%
(i=1,...,n)  \label{Yi}
\end{equation}%
satisfying the eigenvalue equations 
\begin{equation}
Y_{i}E_{\eta }\left( z;q,t\right) =\overline{\eta }_{i}E_{\eta }\left(
z;q,t\right) .  \label{eigen}
\end{equation}%
In (\ref{Yi}) $T_{i}$ denotes the Demazure-Lustig operator, 
\begin{equation}
T_{i}:=t+\frac{tz_{i}-z_{i+1}}{z_{i}-z_{i+1}}\left( s_{i}-1\right) ,
\label{Ti}
\end{equation}%
while%
\begin{equation*}
\omega :=s_{n-1}...s_{1}\tau _{1},
\end{equation*}%
where $s_{i}$ is a transposition operator with the action on functions%
\begin{equation}
\left( s_{i}f\right) \left( z_{1},...,z_{i},z_{i+1},...,z_{n}\right)
:=f\left( z_{1},...,z_{i+1},z_{i},...,z_{n}\right) .  \label{switching}
\end{equation}%
The operator $\tau _{i}$ has the action on functions 
\begin{equation*}
\left( \tau _{i}f\right) \left( z_{1},...,z_{n}\right) :=f\left(
z_{1},...,qz_{i},...,z_{n}\right)
\end{equation*}%
and so corresponds to a $q-$shift of the variable $z_{i}.$ The eigenvalue $%
\overline{\eta }_{i}$ in $\left( \ref{eigen}\right) $ is given by 
\begin{equation}
\overline{\eta }_{i}:=q^{\eta _{i}}t^{-l_{\eta }^{\prime }\left( i\right) },
\label{eigenvalue}
\end{equation}%
where%
\begin{equation}
l_{\eta }^{\prime }(i):=\#\left\{ j<i;\eta _{j}\geq \eta _{i}\right\}
-\#\left\{ j>i;\eta _{j}>\eta _{i}\right\} .  \label{leg colength}
\end{equation}%
Nonsymmetric Macdonald polynomials are of the triangular form 
\begin{equation}
E_{\eta }\left( z;q,t\right) :=z^{\eta }+\sum\limits_{\nu \prec \eta
}b_{\eta \nu }z^{\nu },  \label{form}
\end{equation}%
for coefficients $b_{\eta \nu }$. The notation $z^{\eta }$ denotes the
monomial 
\begin{equation*}
z^{\eta }:=z_{1}^{\eta _{1}}...z_{n}^{\eta _{n}}.
\end{equation*}%
In $\left( \ref{form}\right) $ the coefficient of $z^{\eta }$ has been
chosen to be unity as a normalisation. The ordering $\prec $ is a partial
ordering on compositions having the same modulus, where $\left\vert \eta
\right\vert :=\Sigma _{i=1}^{n}\eta _{i}$ denotes the modulus of $\eta $.
The partial ordering is defined by%
\begin{equation*}
\mu \prec \eta \text{ iff }\mu ^{+}<\eta ^{+}\text{ or in the case }\mu
^{+}=\eta ^{+},\text{ }\mu <\eta
\end{equation*}%
where $\eta ^{+}$ is the unique partition obtained by permuting the
components of $\eta $ and $\mu <\eta $ iff $\mu \neq \eta $ and $\Sigma
_{i=1}^{p}\left( \eta _{i}-\mu _{i}\right) \geq 0$ for all $1\leq p\leq n.$%
\begin{equation*}
\end{equation*}%
Nonsymmetric Macdonald polynomials were first introduced in 1994 [\ref%
{affine},\ref{cherednik}], six years after Macdonald's paper [\ref{macdonald
again}] introducing what are now referred to as symmetric Macdonald
polynomials $P_{\kappa }\left( z;q,t\right) .$ The symmetric Macdonald
polynomials are indexed by partitions $\kappa $ rather than compositions.
The nonsymmetric Macdonald polynomials can be regarded as building blocks of
their symmetric counterparts, as symmetrisation of $E_{\eta }$ gives $%
P_{\eta ^{+}}$. The required symmetrisation operation is defined by%
\begin{equation*}
U^{+}:=\sum\limits_{\sigma \in S_{n}}T_{\sigma },
\end{equation*}%
where $S_{n}$ denotes the set of all permutations of $%
%TCIMACRO{\U{2115} }%
%BeginExpansion
\mathbb{N}
%EndExpansion
^{n}$ and with $\sigma :=s_{i_{l}}...s_{i_{1}}$ the operator $T_{\sigma }$
is specified by%
\begin{equation*}
T_{\sigma }:=T_{i_{l}}...T_{i_{1}},
\end{equation*}%
where $T_{i}$ is defined by $\left( \ref{Ti}\right) .$ The symmetrising
operator allows many fundamental properties of the symmetric Macdonald
polynomials to be deduced as corollaries of the corresponding properties of
the nonsymmetric Macdonald polynomials [\ref{marshall macdonald}]. However,
the converse does not apply, as some special properties of symmetric
Macdonald polynomials have no known nonsymmetric analogues. For example, the
Pieri-type formula $[\ref{macdonald},$ Section VI. 6$]$%
\begin{equation}
e_{r}\left( z\right) P_{\kappa }\left( z;q,t\right) =\sum\limits_{\lambda
}\psi _{\lambda /\kappa }P_{\lambda }\left( z;q,t\right)  \label{pieri}
\end{equation}%
giving the explicit form of the branching coefficients $\psi _{\lambda
/\kappa }$ for the product of $P_{\kappa }\left( z;q,t\right) $ with the $%
r^{th}$ elementary symmetric function 
\begin{equation*}
e_{r}\left( z\right) =\sum\limits_{1\leq i_{1}<...<i_{r}\leq
n}z_{i_{1}}...z_{i_{r}}
\end{equation*}%
has no known non-symmetric analogue. In $\left( \ref{pieri}\right) $ the sum
is over $\lambda $ such that $\lambda /\kappa $ is a vertical $m-$strip.%
\begin{equation*}
\end{equation*}%
Pieri-type formulas themselves have found recent applications in studies of
certain vanishing properties of Macdonald polynomials at $t^{k+1}q^{r-1}=1$ [%
\ref{japanese}]. Furthermore, the dual of $\left( \ref{pieri}\right) $ has
found application in the study of certain probabilistic models related to
the Robinson-Schensted-Knuth corresponence [\ref{last peter}]. 
\begin{equation*}
\end{equation*}%
It is the main objective of this paper to provide the explicit branching
coefficients for the products $z_{i}E_{\eta }\left( z;q,t\right) ,$ $%
e_{1}\left( z\right) E_{\eta }\left( z;q,t\right) $ and $e_{n-1}\left(
z\right) E_{\eta }\left( z;q,t\right) $ in terms of higher order
nonsymmetric polynomials$.$ The latter expansions, which in fact will be
derived as corollaries of the first, are the nonsymmetric analogues of $%
\left( \ref{pieri}\right) $ for $r=1$ and $r=n-1.$ 
\begin{equation*}
\end{equation*}%
That such branching formulas can be derived is suggested by Jack polynomial
theory. Jack polynomials $E_{\eta }\left( z;\alpha \right) $ are the limit $%
q=t^{\alpha },$ $q\rightarrow 1$ of Macdonald polynomials. Marshall [\ref%
{marshall jack}] derived the branching coefficients for the products $%
z_{i}E_{\eta }\left( z;\alpha \right) $ and $e_{1}\left( z\right) E_{\eta
}\left( z;\alpha \right) $ following a strategy of Knop and Sahi [\ref{knop
and sahi}], which proceeds by exploiting the theory of interpolation Jack
polynomials. 
\begin{equation*}
\end{equation*}%
Similarly to the Jack case, the interpolation polynomials play a key role in
deriving the branching coefficients for the product $z_{i}E_{\eta }\left(
z;q,t\right) .$ The nonsymmetric interpolation Macdonald polynomials are
denoted by $E_{\eta }^{\ast }\left( z;q,t\right) $ and can be defined, up to
normalisation, as the unique polynomial of degree $\leq \left\vert \eta
\right\vert $ satisfying 
\begin{equation}
E_{\eta }^{\ast }\left( \overline{\mu }\right) =0,\text{ \ }\left\vert \mu
\right\vert \leq \left\vert \eta \right\vert ,\mu \neq \eta
\label{interpolation}
\end{equation}%
and $E_{\eta }^{\ast }\left( \overline{\eta }\right) \not=0,$ where $%
\overline{\eta }=\left( \overline{\eta }_{1},...,\overline{\eta }_{n}\right) 
$ with $\overline{\eta }_{j}$ specified by (\ref{eigenvalue}). The
interpolation Macdonald polynomials have a triangular expansion in terms of
Macdonald polynomials 
\begin{equation*}
E_{\eta }^{\ast }\left( z;q,t\right) =E_{\eta }\left( z;q,t\right)
+\sum\limits_{\substack{ \left\vert \mu \right\vert \leq \left\vert \eta
\right\vert  \\ \mu \neq \eta }}b_{\eta \mu }E_{\mu }\left( z;q,t\right) ,
\end{equation*}%
for coefficients $b_{\eta \mu }.$ Again, the leading coefficient has been
chosen to be unity as a normalisation. 
\begin{equation*}
\end{equation*}%
The overall strategy for finding the coefficients is to introduce a mapping $%
\Psi $ between $E_{\eta }$ and $E_{\eta }^{\ast }$ that can be used to
intertwine the actions of multiplication by $z_{i}$ on $E_{\eta }$ and a
certain operator $Z_{i}$ on $E_{\eta }^{\ast }.$ Hence, by first determining
an explicit form for the coefficients $c_{\lambda \eta }^{\left\{ i\right\}
} $ in the expansion 
\begin{equation}
Z_{i}E_{\eta }^{\ast }\left( z\right) =\sum\limits_{\nu }c_{\lambda \eta
}^{\left\{ i\right\} }E_{\lambda }^{\ast }\left( z\right)
\label{explicit formula}
\end{equation}%
we can apply the mapping $\Psi $ to obtain an explicit form of the
coefficients of $z_{i}E_{\eta }$ in terms of the $E_{\lambda }$ $\left( 
\text{Section \ref{decomposition}}\right) .$ Using this result we can derive
the explicit formula for the expansion of $e_{1}\left( z\right) E_{\eta }$ $%
\left( \text{Section \ref{pieri one}}\right) .$ The expansion of $%
e_{n-1}\left( z\right) E_{\eta }$ $\left( \text{Section \ref{pieri last}}%
\right) $ then follows from this using the identity $E_{\eta }\left(
z^{-1};q,t\right) =E_{-\eta ^{R}}\left( z;q,t\right) $ [\ref{marshall
macdonald}], where $\eta ^{R}:=\left( \eta _{n},...,\eta _{1}\right) .$ 
\begin{equation*}
\end{equation*}%
The branching coefficients for $z_{i}E_{\eta },$ $e_{1}\left( z\right)
E_{\eta }$ and $e_{n-1}\left( z\right) E_{\eta }$ are given in Propositions %
\ref{final proposition}, \ref{product with elementary function} and \ref%
{en-1}, respectively. As a consequence of finding $e_{1}\left( z\right)
E_{\eta }$ we are able to give an evaluation of the generalised binomial
coefficients $\binom{\eta }{\nu }$ associated with the nonsymmetric
Macdonald polynomials for $\left\vert \eta \right\vert =\left\vert \nu
\right\vert +1$ $\left( \text{Section \ref{product with elementary function}}%
\right) .$ This is given in Proposition \ref{binomial formula}.%
\begin{equation*}
\end{equation*}%
In the final section we take the limit $t=q^{1/\alpha },$ $q\rightarrow 1$
of our result for $e_{1}\left( z\right) E_{\eta }\left(
z;q^{-1},t^{-1}\right) $ to reclaim the known expansion of $E_{\eta }\left(
z;\alpha \right) $ in the theory of nonsymmetric Jack polynomials [\ref%
{marshall jack}].
\begin{equation*}
\end{equation*}%
Note added: After completing this work, and posting it on the arXiv, correspondence was received from Ole Warnaar, pointing out a recent manuscript of Lascoux [\ref{lascoux}], available only on his website, containing results equivalent to our Propositions 7 and 8.

\section{Hecke Operators and the Intertwining Formula\label{preliminaries}}

Hecke operators play an important role in interpolation Macdonald polynomial
theory. They are realisations of the type-A Hecke algebra 
\begin{eqnarray}
\left( H_{i}+1\right) \left( H_{i}-t\right) &=&0  \notag \\
H_{i}H_{i+1}H_{i} &=&H_{i+1}H_{i}H_{i+1},\text{ }i=2,...,n-2
\label{hecke algebra} \\
H_{i}H_{j} &=&H_{j}H_{i},\text{ }\left\vert i-j\right\vert >1.\text{ } 
\notag
\end{eqnarray}%
The Hecke operators of interest, $H_{i},$ are defined by%
\begin{equation}
H_{i}:=\frac{\left( t-1\right) z_{i}}{z_{i}-z_{i+1}}+\frac{z_{i}-tz_{i+1}}{%
z_{i}-z_{i+1}}s_{i},  \label{1a}
\end{equation}%
where $s_{i}$ is specified by $\left( \ref{switching}\right) $. These Hecke
operators appear in the eigenoperators of the interpolation Macdonald
polynomials. The eigenoperators, which mutually commute are defined by 
\begin{equation}
\Xi _{i}:=z_{i}^{-1}+z_{i}^{-1}H_{i}...H_{n-1}\Phi H_{1}...H_{i-1},
\label{XiI}
\end{equation}%
where 
\begin{equation}
\Phi :=\left( z_{n}-t^{-n+1}\right) \Delta \text{ }  \label{Phi}
\end{equation}%
and%
\begin{equation*}
\Delta f\left( z_{1},...,z_{n}\right) =f\left( \frac{z_{n}}{q}%
,z_{1},...,z_{n-1}\right) .
\end{equation*}%
Explicitly, the operators $\Xi _{i}$ satisfy 
\begin{equation}
\Xi _{i}E_{\eta }^{\ast }\left( z;q,t\right) =\overline{\eta }%
_{i}^{-1}E_{\eta }^{\ast }\left( z;q,t\right) ,  \label{eigenfunction}
\end{equation}%
where $\overline{\eta }_{i}$ is given by (\ref{eigenvalue}). The algebraic
relations $\left( \ref{hecke algebra}\right) $ are invariant under the
mapping $H_{i}\longmapsto -H_{i}-1+t.$ Hence, the operators $-\overline{H}%
_{i},$ where%
\begin{equation*}
\overline{H}_{i}:=\frac{\left( t-1\right) z_{i+1}}{z_{i}-z_{i+1}}+\frac{%
z_{i}-tz_{i+1}}{z_{i}-z_{i+1}}s_{i}
\end{equation*}%
are also realisations of the type-A Hecke algebra. These operators appear in
the eigenoperator of $E_{\eta }\left( z;q^{-1},t^{-1}\right) $ according to 
\begin{equation*}
\xi _{i}^{-1}:=\overline{H}_{i}...\overline{H}_{n-1}\Delta H_{1}...H_{i-1}.
\end{equation*}%
By observing 
\begin{equation*}
\Xi _{i}=\xi _{i}^{-1}+\text{degree lowering terms,}
\end{equation*}%
Knop [\ref{knop}] showed that the top homogeneous component of any
interpolation Macdonald polynomial $E_{\eta }^{\ast }\left( z;q,t\right) $
is $E_{\eta }\left( z;q^{-1},t^{-1}\right) .$ Hence, we can define an
isomorphism $\Psi $ mapping each Macdonald polynomial $E_{\eta }\left(
z;q^{-1},t^{-1}\right) $ to its corresponding interpolation polynomial $%
E_{\eta }^{\ast }\left( z;q,t\right) ,$ 
\begin{equation}
\Psi E_{\eta }\left( z;q^{-1},t^{-1}\right) =E_{\eta }^{\ast }\left(
z;q,t\right) .  \label{2a}
\end{equation}%
From this isomorphism we are able to define the important intertwining
formula, Eqn $\left( \ref{inversion}\right) $ below. This is due to Knop [%
\ref{knop}], however in the following an alternative proof is given.

\begin{proposition}
\label{ZiPhi}[\ref{knop}, Theorem 5.1] Define 
\begin{equation}
Z_{i}:=t^{^{- \binom{n}{2}} }\left( z_{i}\Xi _{i}-1\right) \Xi
_{1}...\widehat{\Xi }_{i}...\Xi _{n},  \label{Zi}
\end{equation}%
where the hat superscript on $\widehat{\Xi }_{i}$ denotes the absence of $%
\Xi _{i}$ in the product of operators $\Pi _{j=1}^{n}\Xi _{j}$, and let $M$
be the operator which acts on the subspace of homogeneous polynomials of
degree $d$ by multiplication with $q^{- \binom{d}{2} }.$ With $%
\Psi $ as defined in $\left( \ref{2a}\right) $ we have 
\begin{equation}
Z_{i}\Psi M=\Psi Mz_{i}.  \label{inversion}
\end{equation}
\end{proposition}

\begin{proof}
First consider the action of $Z_{i}$ on $E_{\eta }^{\ast }\left(
z;q,t\right) .$ By the definition of $Z_{i}$ and commutativity of the $\Xi
_{i}$ we have%
\begin{equation}
Z_{i}E_{\eta }^{\ast }\left( z;q,t\right) =(z_{i}-\Xi _{i}^{-1})t^{^{- 
\binom{n}{2} }}\Xi _{1}...\Xi _{n}E_{\eta }^{\ast }\left(
z;q,t\right) .  \label{initial}
\end{equation}%
Using $\left( \ref{eigenfunction}\right) ,$ $\left( \ref{eigenvalue}\right)
, $ then the identity $\Sigma _{i}l_{\eta }^{\prime }\left( i\right) = 
\binom{n}{2}$ we can simplify $\left( \ref{initial}\right) $ to 
\begin{equation}
q^{-\left\vert \eta \right\vert }(z_{i}-\overline{\eta }_{i})E_{\eta }^{\ast
}\left( z;q,t\right) .  \label{final}
\end{equation}%
Since $\left( \ref{final}\right) $ vanishes for all $z=\overline{\lambda },$ 
$\left\vert \lambda \right\vert \leq \left\vert \eta \right\vert ,$ due to $%
\left( \ref{interpolation}\right) ,$ and has degree $\left\vert \eta
\right\vert +1$ we must have 
\begin{equation}
Z_{i}E_{\eta }^{\ast }\left( z;q,t\right) =q^{-\left\vert \eta \right\vert
}\sum\limits_{\lambda :\left\vert \lambda \right\vert =\left\vert \eta
\right\vert +1}c_{\lambda \eta }^{\left\{ i\right\} }E_{\lambda }^{\ast
}\left( z;q,t\right) ,  \label{ZiEn}
\end{equation}%
for some coefficients $c_{\lambda \eta }^{\left\{ i\right\} }.$ Equating the
leading terms of (\ref{final}) and the right hand side of $\left( \ref{ZiEn}%
\right) $ gives%
\begin{equation}
z_{i}E_{\eta }\left( z;q^{-1},t^{-1}\right) =\sum\limits_{\lambda
:\left\vert \lambda \right\vert =\left\vert \eta \right\vert +1}c_{\lambda
\eta }^{\left\{ i\right\} }E_{\lambda }\left( z;q^{-1},t^{-1}\right) .
\label{3a}
\end{equation}%
Applying the action of $\Psi M$ to both sides of $\left( \ref{3a}\right) $
and using $\left( \ref{2a}\right) $ shows 
\begin{equation}
\Psi Mz_{i}E_{\eta }\left( z;q^{-1},t^{-1}\right) =q^{-\tbinom{\left\vert
\eta \right\vert +1}{2}}\sum\limits_{\lambda :\left\vert \lambda
\right\vert =\left\vert \eta \right\vert +1}c_{\lambda \eta }^{\left\{
i\right\} }E_{\lambda }^{\ast }\left( z;q,t\right) .  \label{next one}
\end{equation}%
Using $\left( \ref{ZiEn}\right) ,$ the right hand side of $\left( \ref{next
one}\right) $ can be simplified to 
\begin{equation}
q^{-\tbinom{\left\vert \eta \right\vert }{2}}Z_{i}E_{\eta }^{\ast }\left(
z;q,t\right) .  \label{ZiEnStar}
\end{equation}%
By recalling the action of $M$ and again using $\left( \ref{2a}\right) $ we
obtain%
\begin{equation}
\Psi Mz_{i}E_{\eta }\left( z;q^{-1},t^{-1}\right) =Z_{i}\Psi ME_{\eta
}\left( z;q^{-1},t^{-1}\right) .  \label{intertwining result a}
\end{equation}%
Finally, since the $\left\{ E_{\eta }\right\} $ form a basis for analytic
functions in $\left\{ z^{\eta }\right\} $ it follows that the intertwining
property $\left( \ref{inversion}\right) $ holds generally.
\end{proof}

\begin{corollary}
\label{corollary 3}We have%
\begin{equation}
z_{i}E_{\eta }\left( z;q^{-1},t^{-1}\right) =q^{\left\vert \eta \right\vert
}\Psi ^{-1}Z_{i}E_{\eta }^{\ast }\left( z;q,t\right) .  \label{the result}
\end{equation}
\end{corollary}

\begin{proof}
Follows from $\left( \ref{next one}\right) $ and $\left( \ref{ZiEnStar}%
\right) .$
\end{proof}

\section{The Product $z_{i}E_{\protect\eta }$\label{decomposition}}

The previous corollary indicates that the next step towards finding the
decomposition of $z_{i}E_{\eta }$ is to determine an explicit formula for $%
Z_{i}E_{\eta }^{\ast }.$ The latter can be deduced as a corollary of the
following lemma, specifying the expansion of $\left( z_{i}\Xi _{i}-1\right)
f\left( z\right) ,$ where according to $\left( \ref{XiI}\right) $ $z_{i}\Xi
_{i}-1:=H_{i}...H_{n-1}\Phi H_{1}...H_{i-1}.$

\begin{lemma}
\label{lemma 1}Let $\widetilde{Z}_{i}=H_{i}...H_{n-1}\Phi H_{1}...H_{i-1}.$
The action of $\widetilde{Z}_{i}$ on $f\left( z\right) $ is given by%
\begin{equation}
\widetilde{Z}_{i}f\left( z\right) =\sum\limits_{\substack{ I\subseteq
\left\{ 1,...,n\right\}  \\ i\in I}}r_{I}^{\left\{ i\right\} }\left(
z\right) f\left( Iz\right) .  \label{4a}
\end{equation}%
Here the rational function $r_{I}^{\left\{ i\right\} }\left( z\right)$ can be expressed as 
\begin{equation}
r_{I}^{\left\{ i\right\} }\left( z\right) =\chi _{I}^{\left\{ i\right\}
}\left( z\right) A_{I}\left( z\right) B_{I}\left( z\right)  \label{4b}
\end{equation}%
where 
\begin{eqnarray}
I&=&\left\{ t_{1},...,t_{s}\right\} ,1\leq t_{1}<...t_{s}\leq n,  \label{I} \\
A_{I}\left( z\right) &=&\widehat{a}\left( \frac{z_{t_{s}}}{q}%
,z_{t_{1}}\right) \prod\limits_{u=1}^{s-1}\widehat{a}\left(
z_{t_{u}},z_{t_{u+1}}\right)  \label{Ai} \\
B_{I}\left( z\right) &=&\left( z_{t_{s}}-t^{-n+1}\right) \left(
\prod\limits_{j=1}^{t_{1}-1}\widehat{b}\left( \frac{z_{t_{s}}}{q}%
,z_{j}\right) \right)  \notag \\
&&\times \left( \prod\limits_{u=1}^{s}\prod\limits_{j=t_{u}+1}^{t_{u+1}-1}%
\widehat{b}\left( z_{t_{u}},z_{j}\right) \right) ,\text{ }t_{s+1}:=n+1
\label{Bi} \\
\chi _{I}^{\left\{ i\right\} }\left( z\right) &=&\left\{ 
\begin{tabular}{ll}
$\frac{1}{\widehat{a}\left( z_{t_{k-1}},z_{i}\right) }$ & $;i=t_{k},$ $%
k=2,...,s$ \\ 
$\frac{1}{\widehat{a}\left( \frac{z_{t_{s}}}{q},z_{i}\right) }$ & $;i=t_{1},$%
\end{tabular}%
,\right.  \label{Xi} \\
\end{eqnarray}%
and $Iz$ is defined as 
\begin{equation*}
\left( Iz\right) _{i}=\left\{ 
\begin{tabular}{ll}
$z_{t_{u-1}}$ & $;i=t_{u},$ if $u=2,...,s$ \\ 
$\frac{z_{t_{s}}}{q}$ & $;i=t_{1}$ \\ 
$z_{i}$ & $;i\notin I.$%
\end{tabular}%
\right. .
\end{equation*}%
The quantities $\widehat{a}\left( x,y\right) $ and $\widehat{b}\left(
x,y\right) $ are defined in $\left( \ref{a and b}\right) $ below.
\end{lemma}

\begin{proof}
Using ($\ref{1a})$ the action of $H_{i}$ on $f\left( z\right) $ can be
expressed as%
\begin{equation*}
H_{i}f\left( z\right) =\widehat{a}\left( x,y\right) f\left( z\right) +%
\widehat{b}\left( x,y\right) s_{i}f\left( z\right) ,
\end{equation*}%
where%
\begin{equation}
\widehat{a}\left( x,y\right) :=\frac{\left( t-1\right) x}{x-y},\text{ }%
\widehat{b}\left( x,y\right) :=\frac{x-ty}{x-y}.  \label{a and b}
\end{equation}%
Hence $\widetilde{Z}_{i}$ can be written as 
\begin{eqnarray}
&&\left( \widehat{a}\left( z_{i},z_{i+1}\right) +\widehat{b}\left(
z_{i},z_{i+1}\right) s_{i}\right) ...\left( \widehat{a}\left(
z_{n-1},z_{n}\right) +\widehat{b}\left( z_{n-1},z_{n}\right) s_{n-1}\right) 
\notag \\
&&\qquad \times \Phi \left( \widehat{a}\left( z_{1},z_{2}\right) +\widehat{b}\left(
z_{1},z_{2}\right) s_{1}\right) ...\left( \widehat{a}\left(
z_{i-1},z_{i}\right) +\widehat{b}\left( z_{i-1},z_{i}\right) s_{i-1}\right) .
\label{product}
\end{eqnarray}%
Let 
\begin{equation*}
K_{I}^{\left\{ i\right\} }=s_{i}...\widehat{s}_{t_{r+1}-1}...\widehat{s}%
_{t_{s}-1}...s_{n-1}\Delta s_{1}...\widehat{s}_{t_{1}}...\widehat{s}%
_{t_{r-1}}...s_{i-1},\text{ for }i\in I,
\end{equation*}%
where $1\leq t_{1}<...<t_{r}=i<t_{r+1}<...<t_{s}\leq n$, the hat superscript
used as in Section $\ref{preliminaries}$ to denote the absence of the
corresponding operators and $I$ as defined in the statement of the result.
It is clear that the expansion of $\widetilde{Z}_{i}$ will be of the form 
\begin{equation*}
\widetilde{Z}_{i}=\sum\limits_{\substack{ I\subseteq \left\{
1,...,n\right\}  \\ i\in I}}r_{I}^{\left\{ i\right\} }\left( z\right)
K_{I}^{\left\{ i\right\} }
\end{equation*}%
for coefficients $r_{I}^{\left\{ i\right\} }\left( z\right) $ involving $%
\widehat{a}\left( x,y\right) $ and $\widehat{b}\left( x,y\right) .$ Further,
it is easily verified that $K_{I}^{\left\{ i\right\} }f\left( z\right)
=f\left( Iz\right) .$ The coefficients $r_{I}^{\left\{ i\right\} }\left(
z\right) $ are found by considering the individual terms in the expansion of 
$\left( \ref{product}\right) .$ Due to the need to commute the transposition
operators $s_{i}$ through to the right the final formula is more simply
obtained by expanding $\left( \ref{product}\right) $ termwise from the
right. Inevitably, the exercise is rather tedious, however it can be
structured somewhat by considering four disjoint classes of sets $I$%
\begin{eqnarray*}
I_{1} &=&\left\{ i\right\} , \\
I_{2} &=&\left\{ ...,i\right\} , \\
I_{3} &=&\left\{ i,...\right\} , \\
I_{4} &=&\left\{ ...,i,...\right\} ,
\end{eqnarray*}%
which exhaust all possibilities. This cataloguing allows the coefficients of
the corresponding four forms of $K_{I}^{\left\{ i\right\} }$ to be
considered separately and the result is more easily observed. Explicitly,
the four forms of $K_{I}^{\left\{ i\right\} }$ are 
\begin{eqnarray*}
&&K_{I_{1}}^{\left\{ i\right\} }=s_{i}...s_{n-1}\Delta s_{1}...s_{i-1}, \\
&&K_{I_{2}}^{\left\{ i\right\} }=\text{ }s_{i}...s_{n-1}\Delta s_{1}...%
\widehat{s}_{t_{1}}...\widehat{s}_{t_{r-1}}...s_{i-1}, \\
&&K_{I_{3}}^{\left\{ i\right\} }=\text{ }s_{i}...\widehat{s}_{t_{r+1}-1}...%
\widehat{s}_{t_{s}-1}...s_{n-1}\Delta s_{1}...s_{i-1}, \\
&&K_{I_{4}}^{\left\{ i\right\} }=\text{ }s_{i}...\widehat{s}_{t_{r+1}-1}...%
\widehat{s}_{t_{s}-1}...s_{n-1}\Delta s_{1}...\widehat{s}_{t_{1}}...\widehat{%
s}_{t_{r-1}}...s_{i-1}.
\end{eqnarray*}%
In relation to $K_{I_{2}}$, $K_{I_{4}}$ the coefficient of $s_{1}...\widehat{%
s}_{t_{1}}...\widehat{s}_{t_{r-1}}...s_{i-1}$ in the partial expansion, that
is terms to the right of $\Phi ,$ of $\left( \ref{product}\right) $ is%
\begin{eqnarray*}
&&\widehat{a}\left( z_{1},z_{t_{1}+1}\right) \prod\limits_{u=1}^{r-2}%
\widehat{a}\left( z_{t_{u}+1},z_{t_{u+1}+1}\right) \\
&&\qquad \times \prod\limits_{j=1}^{t_{1}-1}\widehat{b}\left( z_{1},z_{j+1}\right)
\prod\limits_{u=1}^{r-1}\prod\limits_{j=t_{u}+1}^{t_{u+1}-1}\widehat{b}%
\left( z_{t_{u}+1},z_{j+1}\right) .
\end{eqnarray*}%
Hence the coefficient of $\Delta s_{1}...\widehat{s}_{t_{1}}...\widehat{s}%
_{t_{r-1}}...s_{i-1}$ will be%
\begin{eqnarray*}
&&\left( z_{n}-t^{-n+1}\right)\widehat{a}\left( \frac{z_{n}}{q}%
,z_{t_{1}}\right) \prod\limits_{u=1}^{r-2}\widehat{a}\left(
z_{t_{u}},z_{t_{u+1}}\right) \\
&& \qquad \times \prod\limits_{j=1}^{t_{1}-1}\widehat{b}\left( \frac{z_{n}}{q}%
,z_{j}\right) \prod\limits_{u=1}^{r-1}\prod\limits_{j=t_{u}+1}^{t_{u+1}-1}%
\widehat{b}\left( z_{t_{u}},z_{j}\right) .
\end{eqnarray*}%
Similarly, for $K_{I_{1}}$, $K_{I_{3}}$, the coefficient of $s_{1}...s_{i-1}$
and $\Delta s_{1}...s_{i-1}$ are%
\begin{equation*}
\prod\limits_{j=1}^{i-1}\widehat{b}\left( z_{1},z_{j+1}\right)
\end{equation*}%
and%
\begin{equation*}
\left( z_{n}-t^{-n+1}\right) \prod\limits_{j=1}^{i-1}\widehat{b}\left( 
\frac{z_{n}}{q},z_{j}\right) ,
\end{equation*}%
respectively. The final $r_{I_{j}}\left( z\right) ^{\prime }s$ are found by
continuing the expansion of $\left( \ref{product}\right) $ from the right
and considering the four forms of $K_{I}$ separately. Thus we find that 
\begin{eqnarray*}
r_{I_{1}}^{\left\{ i\right\} }\left( z\right) &=&\frac{A_{I}\left( z\right)
B_{I}\left( z\right) }{\widehat{a}\left( \frac{z_{i}}{q},z_{i}\right) } \\
r_{I_{2}}^{\left\{ i\right\} }\left( z\right) &=&\frac{A_{I}\left( z\right)
B_{I}\left( z\right) }{\widehat{a}\left( z_{t_{r-1}},z_{t_{r}}\right) } \\
r_{I_{3}}^{\left\{ i\right\} }\left( z\right) &=&\frac{A_{I}\left( z\right)
B_{I}\left( z\right) }{\widehat{a}\left( \frac{z_{t_{s}}}{q}%
,z_{t_{1}}\right) } \\
r_{I_{4}}^{\left\{ i\right\} }\left( z\right) &=&\frac{A_{I}\left( z\right)
B_{I}\left( z\right) }{\widehat{a}\left( z_{t_{r-1}},z_{t_{r}}\right) },
\end{eqnarray*}%
where $A_{I}\left( z\right) $ and $B_{I}\left( z\right) $ are defined by $%
\left( \ref{Ai}\right) $ and $\left( \ref{Bi}\right) ,$ respectively$.$
After recalling the definition of $\chi _{I}^{\left\{ i\right\} }$ given
above, the sought explicit formula $\left( \ref{4b}\right) $ follows.
\end{proof}

\begin{corollary}
\label{corollary}We have 
\begin{equation}
Z_{i}E_{\eta }^{\ast }\left( z\right) =q^{-\left\vert \eta \right\vert }%
\overline{\eta }_{i}\sum\limits_{\substack{ I\subseteq \left\{
1,...,n\right\}  \\ i\in I}}r_{I}^{\left\{ i\right\} }\left( z\right)
E_{\eta }^{\ast }\left( Iz\right) .  \label{corollary formula}
\end{equation}
\end{corollary}

\begin{proof}
Follows after recalling from $\left( \ref{Zi}\right) $ that 
\begin{equation*}
Z_{i}:=t^{^{- \binom{n}{2} }}\left( z_{i}\Xi _{i}-1\right) \Xi
_{1}...\widehat{\Xi }_{i}...\Xi _{n}.
\end{equation*}
\end{proof}

Together Proposition \ref{ZiPhi} and Corollary \ref{corollary} allow us to
derive an initial expansion $z_{i}E_{\eta }\left( z;q^{-1},t^{-1}\right) $
in terms of the Macdonald polynomials of degree $\left\vert \eta \right\vert
+1$.

\begin{proposition}
We have%
\begin{equation}
z_{i}E_{\eta }\left( z;q^{-1},t^{-1}\right) =\overline{\eta }%
_{i}q^{-\left\vert \eta \right\vert }\sum\limits_{\left\vert \lambda
\right\vert =\left\vert \eta \right\vert +1}\sum\limits_{_{\substack{ %
I\subseteq \left\{ 1,...,n\right\}  \\ i\in I}}}\frac{r_{I}^{\left\{
i\right\} }\left( \overline{\lambda }\right) E_{\eta }^{\ast }\left( I%
\overline{\lambda }\right) }{E_{\lambda }^{\ast }\left( \overline{\lambda }%
\right) }E_{\lambda }\left( z;q^{-1},t^{-1}\right) .
\label{initial expansion}
\end{equation}
\end{proposition}

\begin{proof}
By the vanishing properties of $E_{\eta }^{\ast },$ $\left( \ref%
{interpolation}\right) ,$ when the right hand sides of $\left( \ref{ZiEn}%
\right) $ and $\left( \ref{corollary formula}\right) $ equated and evaluated
at $z=\overline{\lambda }$ we obtain%
\begin{equation*}
c_{\lambda \eta }^{\left\{ i\right\} }=\overline{\eta }_{i}\sum\limits_{ 
_{\substack{ I\subseteq \left\{ 1,...,n\right\}  \\ i\in I}}}\frac{%
r_{I}^{\left\{ i\right\} }\left( \overline{\lambda }\right) E_{\eta }^{\ast
}\left( I\overline{\lambda }\right) }{E_{\lambda }^{\ast }\left( \overline{%
\lambda }\right) }.
\end{equation*}%
Substituting this back into $\left( \ref{ZiEn}\right) $ and applying
Corollary \ref{corollary 3} gives $\left( \ref{initial expansion}\right) .$
\end{proof}

\begin{equation*}
\end{equation*}

The formula $\left( \ref{initial expansion}\right) $ can be improved by
three simplifications. The first is to restrict the summation in $\left( \ref%
{initial expansion}\right) $ by removing a number of vanishing terms. For
this we require the following two propositions, and associated definitions.

\begin{proposition}
\label{comaximal}Let $I=\left\{ t_{1},...,t_{s}\right\} $ with $1\leq
t_{1}<...<t_{s}\leq n$ and $I\neq \emptyset .$ We call $I$ comaximal with
respect to $\lambda $ iff:%
\begin{equation*}
\begin{tabular}{lll}
$\left( 1\right) $ & $\lambda _{j}\neq \lambda _{t_{u}},$ & $%
 j=t_{u}+1,...,t_{u+1}-1,\left( \text{ }u=1,...,s;\text{ }t_{s+1}=n+1\right)
; $ \\ 
$\left( 2\right) $ & $\lambda _{j}\neq \lambda _{t_{s}}-1,$ & $\text{ }%
j=1,...,t_{1}-1;$ \\ 
$\left( 3\right) $ & $\lambda _{t_{s}}\neq 0.$ & 
\end{tabular}%
\end{equation*}%
If $I$ is not comaximal with respect to $\lambda $ then $r_{I}^{\left\{
i\right\} }\left( \overline{\lambda }\right) =0.$
\end{proposition}

\begin{proof}
If any one of the three conditions in the definition of $I$ comaximal with
respect to $\lambda $ fail, then $\widehat{B}_{I}\left( \overline{\lambda }%
\right) =0$ and therefore $r_{I}^{\left\{ i\right\} }\left( \overline{%
\lambda }\right) =0.$
\end{proof}

\begin{proposition}
\label{maximal}Let $I=\left\{ t_{1},...,t_{s}\right\} $ with $1\leq
t_{1}<...<t_{s}\leq n$ and $I\neq \emptyset .$ We call $I$ maximal with
respect to $\lambda $ iff%
\begin{equation}
\begin{tabular}{lll}
$\left( 1\right) $ & $\lambda _{j}\neq \lambda _{t_{u}},$ & $\text{ }%
j=t_{u-1}+1,...,t_{u}-1$ $\left( u=1,...,s;\text{ }t_{0}:=0\right) ;$ \\ 
$\left( 2\right) $ & $\lambda _{j}\neq \lambda _{t_{1}}+1,$ & $\text{ }%
j=t_{s}+1,...,n.$%
\end{tabular}
\label{maximal conditions}
\end{equation}%
Also define the composition $c_{I}\left( \lambda \right) $ for such a set $I$
by%
\begin{equation}
\left( c_{I}\left( \lambda \right) \right) _{j}=\left\{ 
\begin{tabular}{ll}
$\lambda _{t_{k+1}}$ & $;j=t_{k},$ if $k=1,...,s-1$ \\ 
$\lambda _{t_{1}}+1$ & $;j=t_{s}$ \\ 
$\lambda _{j}$ & $;j\notin I.$%
\end{tabular}%
\right. .  \label{c_I}
\end{equation}%
Set $I$ is comaximal with respect to $\lambda $ iff there exists a
composition $\nu $ such that $I$ is maximal with respect to $\nu ,$ $\lambda
=c_{I}\left( \nu \right) $ and $I\overline{\lambda }=\overline{\nu }.$
\end{proposition}

\begin{proof}
Follows from the definitions.
\end{proof}

\begin{equation*}
\end{equation*}%
It is shown in [\ref{knop}] that it is only these maximal subsets which give
distinct compositions $\lambda .$ Thus it is convenient to introduce the set 
$\mathbb{J}_{\eta }^{I}$ of maximal subsets%
\begin{equation}
\mathbb{J}_{\eta }^{I}:=\left\{ I:I\text{ is maximal with respect to }\eta
\right\}  \label{Jn1}
\end{equation}%
and the corresponding set of compositions%
\begin{equation*}
\mathbb{J}_{\eta }^{\lambda }:=\left\{ \lambda :\lambda =c_{I}\left( \eta
\right) ,\text{ }I\in \mathbb{J}_{\eta }^{I}\right\} .
\end{equation*}

\begin{corollary}
\label{corollary 2}If $I$ is comaximal with respect to $\lambda $ then $%
E_{\eta }^{\ast }\left( I\overline{\lambda }\right) \neq 0$ iff $I$ is
maximal with respect to $\eta .$
\end{corollary}

\begin{proof}
Follows from Proposition \ref{maximal} and the vanishing properties of $%
E_{\eta }^{\ast }$ $\left( \ref{interpolation}\right) .$
\end{proof}

\begin{equation*}
\end{equation*}%
Using these results we can begin to simplify $\left( \ref{initial expansion}%
\right) .$

\begin{proposition}
We have 
\begin{equation}
z_{i}E_{\eta }\left( z;q^{-1},t^{-1}\right) =\overline{\eta }%
_{i}\sum\limits _{\substack{ I\in \mathbb{J}_{\eta }^{I},\text{ }i\in I  \\ %
c_{I}\left( \eta \right) =\lambda }}\frac{r_{I}^{\left\{ i\right\} }\left( 
\overline{c_{I}\left( \eta \right) }\right) E_{\eta }^{\ast }\left( 
\overline{\eta }\right) }{E_{c_{I}\left( \eta \right) }^{\ast }\left( 
\overline{c_{I}\left( \eta \right) }\right) }E_{c_{I}\left( \eta \right)
}\left( z;q^{-1},t^{-1}\right) .  \label{second expanstion}
\end{equation}
\end{proposition}

\begin{proof}
Using Proposition \ref{comaximal} we can restrict the second summation of $%
\left( \ref{initial expansion}\right) $ to the sets $I$ that are comaximal
with respect to $\lambda .$ Proposition $\ref{maximal}$ allows us to
restrict the sum further to sets $I$ that are maximal with respect to $\eta $
and hence to $\lambda $ of the form $\lambda =c_{I}\left( \eta \right) ,$
giving the required result.
\end{proof}

\begin{equation*}
\end{equation*}

The second simplification is made by giving an evaluation formula for $%
E_{\eta }^{\ast }\left( \overline{\eta }\right) $. The derivation draws upon
areas of Macdonald polynomial theory not used elsewhere in this work. Hence
to avoid a long deviation from the overall goal, the reader is referred to $%
\left[ \ref{sahi}\right] $ for the details of such results.

\begin{proposition}
\label{evaluation copy(1)}We have%
\begin{equation}
E_{\eta }^{\ast }\left( \overline{\eta }\right) :=k_{\eta }=\left(
\prod\limits_{i=1}^{n}\overline{\eta }_{i}^{\eta _{i}}\right) d_{\eta
}^{\prime }\left( q^{-1},t^{-1}\right) ,  \label{evaluation}
\end{equation}%
where%
\begin{equation}
d_{\eta }^{\prime }\left( q^{-1},t^{-1}\right) :=\prod\limits_{\left(
i,j\right) =s\in \text{diag}\left( \eta \right) }\left( 1-q^{-a_{\eta
}\left( s\right) -1}t^{-l_{\eta }\left( s\right) }\right) ,  \label{d dash}
\end{equation}%
where diag$\left( \eta \right) :=\left\{ \left( i,j\right) \in 
%TCIMACRO{\U{2124} }%
%BeginExpansion
\mathbb{Z}
%EndExpansion
^{2},1\leq j\leq \eta _{i}\right\} .$ The quantities $a_{\eta }\left(
s\right) $ and $l_{\eta }\left( s\right) $ are the arm and leg length
respectively and defined by 
\begin{equation}
a_{\eta }\left( s\right) =\eta _{i}-j\text{ and }l_{\eta }\left( s\right)
=\#\left\{ k>i;j\leq \eta _{k}\leq \eta _{i}\right\} +\#\left\{ k<i;j\leq
\eta _{k}+1\leq \eta _{i}\right\} .  \label{arm and leg}
\end{equation}
\end{proposition}

\begin{proof}
Use will be made of the operations $s_{i}$ and $\Phi $ defined to act on
functions of $n$ variables by $\left( \ref{switching}\right) $ and $\left( %
\ref{Phi}\right) ,$ with their actions now on compositions. The action of $%
s_{i}$ on $\eta $ is to exchange parts in positions $i$ and $i+1,$ while $%
\Phi $ acts on compositions according to 
\begin{equation*}
\Phi \eta :=\left( \eta _{2},...,\eta _{n},\eta _{1}+1\right) .
\end{equation*}%
These operators can generate all compositions recursively, starting with $%
\left( 0,...,0\right) ,$ and allow $\left( \ref{evaluation}\right) $ to be
proved inductively. Clearly, when $\eta =(0,...,0)$ we have $k_{\eta
}=1=E_{\eta }^{\ast }\left( \overline{\eta };q,t\right) ,$ which establishes
the base case. Assume for $\eta $ general $E_{\eta }^{\ast }\left( \overline{%
\eta }\right) =k_{\eta }$. Our task is to deduce from this that%
\begin{equation}
E_{s_{i}\eta }^{\ast }\left( \overline{s_{i}\eta };q,t\right) =k_{s_{i}\eta }
\label{switching1}
\end{equation}%
and%
\begin{equation}
E_{\Phi \eta }^{\ast }\left( \overline{\Phi \eta };q,t\right) =k_{\Phi \eta
}.  \label{raising}
\end{equation}%
To show $\left( \ref{switching1}\right) $ we must consider the cases $\eta
_{i}<\eta _{i+1}$ and $\eta _{i}>\eta _{i+1}$ separately. We begin with the
case $\eta _{i}<\eta _{i+1}.$ To relate $E_{s_{i}\eta }^{\ast }\left( 
\overline{s_{i}\eta }\right) $ to $E_{\eta }^{\ast }\left( \overline{\eta }%
\right) $ we consider two different perspectives on the computation of $%
H_{i}E_{\eta }^{\ast }\left( z\right) $. The first is found by recognising
that $H_{i}=T_{i}^{-1}[t^{-1}],$ where $T_{i}^{-1}$ is the inverse of the
Demazure-Lusztig operator $T_{i}$ defined by $\left( \ref{Ti}\right) .$ From 
$\left( \ref{Ti}\right) $ and the quadratic relation of $\left( \ref{hecke
algebra}\right) $ we have%
\begin{equation}
T_{i}^{-1}:=t^{-1}-1+t^{-1}T_{i}.  \label{ti inv}
\end{equation}%
Taking the known result $[\ref{peter 3}]$%
\begin{equation*}
T_{i}^{-1}E_{\eta }\left( z;q,t\right) =\frac{t^{-1}-1}{1-\delta _{i,\eta
}\left( q,t\right) }E_{\eta }\left( z;q,t\right) +E_{s_{i}\eta }\left(
z;q,t\right) \text{, when }\eta _{i}<\eta _{i+1},
\end{equation*}%
where $\delta _{i,\eta }$ $:=\overline{\eta }_{i}/\overline{\eta }_{i+1},$
and replacing $\left( q,t\right) $ by $\left( q^{-1},t^{-1}\right) $ allows
us to apply the mapping $\Psi $ $\left( \ref{2a}\right) $ to both sides of
the equation. Making use of the fact that $\Psi $ commutes with $H_{i}$ [\ref%
{knop}, Section $5$] we obtain%
\begin{equation}
H_{i}E_{\eta }^{\ast }\left( z;q,t\right) =\frac{t-1}{1-\delta _{i,\eta
}\left( q^{-1},t^{-1}\right) }E_{\eta }^{\ast }\left( z;q,t\right)
+E_{s_{i}\eta }^{\ast }\left( z;q,t\right) .  \label{first hi}
\end{equation}%
\begin{equation*}
\end{equation*}%
The second perspective is obtained directly from definition $\left( \ref{1a}%
\right) $ which gives%
\begin{equation}
H_{i}E_{\eta }^{\ast }\left( z;q,t\right) =\frac{\left( t-1\right) z_{i}}{%
z_{i}-z_{i+1}}E_{\eta }^{\ast }\left( z;q,t\right) +\frac{z_{i}-tz_{i+1}}{%
z_{i}-z_{i+1}}E_{\eta }^{\ast }\left( s_{i}z;q,t\right) .  \label{HiEn again}
\end{equation}%
Equating the right hand sides of $\left( \ref{first hi}\right) $ and $\left( %
\ref{HiEn again}\right) $ and evaluating at $z=\overline{s_{i}\eta }$ we
obtain 
\begin{equation*}
\frac{1-t\delta _{i,\eta }^{-1}\left( q^{-1},t^{-1}\right) }{1-\delta
_{i,\eta }^{-1}\left( q^{-1},t^{-1}\right) }=\frac{E_{s_{i}\eta }^{\ast
}\left( \overline{s_{i}\eta }\right) }{E_{\eta }^{\ast }\left( \overline{%
\eta }\right) }.
\end{equation*}%
Since for $\eta _{i}<\eta _{i+1}$ $[\ref{sahi}]$, 
\begin{equation*}
\frac{d_{s_{i}\eta }^{\prime }\left( q,t\right) }{d_{\eta }^{\prime }\left(
q,t\right) }=\frac{1-\delta _{i,\eta }^{-1}\left( q,t\right) }{%
1-t^{-1}\delta _{i,\eta }^{-1}\left( q,t\right) },
\end{equation*}%
we have 
\begin{equation*}
\frac{k_{s_{i}\eta }}{k_{\eta }}=\frac{d_{s_{i}\eta }^{\prime }\left(
q^{-1},t^{-1}\right) }{d_{\eta }^{\prime }\left( q^{-1},t^{-1}\right) }=%
\frac{E_{s_{i}\eta }^{\ast }\left( \overline{s_{i}\eta }\right) }{E_{\eta
}^{\ast }\left( \overline{\eta }\right) }.
\end{equation*}%
Hence $E_{\eta }^{\ast }\left( \overline{\eta }\right) =k_{\eta }$ implies $%
E_{s_{i}\eta }^{\ast }\left( \overline{s_{i}\eta }\right) =k_{s_{i}\eta }.$
The case where $\eta _{i}>\eta _{i+1}$ is proven similarly. 
\begin{equation*}
\end{equation*}%
The first step to showing $\left( \ref{raising}\right) $ is to consider the
vanishing properties of $\left( \Phi E_{\eta }^{\ast }\right) \left(
z\right) .$ By Knop [\ref{knop}, Corollary $3.3$] if $\left\vert \lambda
\right\vert \leq \left\vert \eta \right\vert $ then $\left( \Phi E_{\eta
}^{\ast }\right) \left( \overline{\lambda }\right) $ is a linear combination
of $E_{\eta }^{\ast }\left( \overline{\nu }\right) $ for $\left\vert \nu
\right\vert \leq \left\vert \lambda \right\vert $ and hence $\left( \Phi
E_{\eta }^{\ast }\right) \left( \overline{\lambda }\right) $ vanishes for $%
\left\vert \lambda \right\vert \leq \left\vert \eta \right\vert $. Now, $%
\left( \Phi E_{\eta }^{\ast }\right) \left( z\right) $ is a polynomial of
degree $\left\vert \eta \right\vert +1$ and of the form 
\begin{equation}
\left( z_{n}-t^{-n+1}\right) E_{\eta }^{\ast }\left( \frac{z_{n}}{q}%
,z_{1},...,z_{n-1}\right) .  \label{polynomial}
\end{equation}%
If $\eta _{n}=0$ \ then $\overline{\eta }_{n}=t^{-n+1}$ and $\left( \ref%
{polynomial}\right) $ is equal to zero. If $\eta _{n}\not=0$ then $\left(
\Phi E_{\eta }^{\ast }\right) \left( \overline{\lambda }\right) $ can be
written as%
\begin{equation*}
\left( \overline{\lambda _{n}}-t^{-n+1}\right) E_{\eta }^{\ast }\left( 
\overline{\nu }\right) ,
\end{equation*}%
where $\nu =\left( \lambda _{n}-1,\lambda _{1},...,\lambda _{n-1}\right) $,
so that $\lambda =\Phi \nu .$ Since $\left\vert \nu \right\vert =\left\vert
\eta \right\vert ,$ by the vanishing properties of $E_{\eta }^{\ast }\left(
z\right) $ if $\lambda \neq \Phi \eta ,$ then $\nu \neq \eta $ and
consequently $\left( \Phi E_{\eta }^{\ast }\right) \left( \overline{\nu }%
\right) =0.$ From these vanishing properties it follows that $\left( \Phi
E_{\eta }^{\ast }\right) \left( z\right) $ is a multiple of $E_{\Phi \eta
}^{\ast }\left( z\right) .$ A computation gives the coefficient of $z^{\Phi
\eta }$ in $(\Phi E_{\eta }^{\ast })(z;q,t)$ to be $q^{-\eta _{1}}$ and so $%
E_{\Phi \eta }^{\ast }\left( z\right) =q^{-\eta _{1}}(\Phi E_{\eta }^{\ast
})(z;q,t).$ By evaluating $\left( \Phi E_{\eta }^{\ast }\right) \left(
z\right) $ at $z=\overline{\Phi \eta }$ and rearranging we obtain 
\begin{equation}
\frac{E_{\Phi \eta }^{\ast }\left( \overline{\Phi \eta }\right) }{E_{\eta
}^{\ast }\left( \overline{\eta }\right) }=q^{-\eta _{1}}\left( q\overline{%
\eta }_{1}-t^{-n+1}\right) .  \label{qt}
\end{equation}%
Now by Sahi $[\ref{sahi}]$ we have 
\begin{equation*}
\frac{d_{\Phi \eta }^{\prime }\left( q,t\right) }{d_{\eta }^{\prime }\left(
q,t\right) }=1-q^{\eta _{1}+1}t^{n-1-l_{\eta }^{\prime }\left( 1\right) }.
\end{equation*}%
Using this and the definition of $k_{\eta }$ we can simplify $\left( \ref{qt}%
\right) $ to 
\begin{equation}
\frac{E_{\Phi \eta }^{\ast }\left( \overline{\Phi \eta }\right) }{E_{\eta
}^{\ast }\left( \overline{\eta }\right) }=q^{2\eta _{1}+1}t^{-l_{\eta
}^{\prime }\left( 1\right) }\frac{d_{\Phi \eta }^{\prime }\left(
q^{-1},t^{-1}\right) }{d_{\eta }^{\prime }\left( q^{-1},t^{-1}\right) }=%
\frac{k_{\Phi \eta }}{k_{\eta }}.  \label{fraction}
\end{equation}%
Where to obtain the final equality the fact that 
\begin{equation*}
\frac{\Phi \left( \prod\limits_{i=1}^{n}\overline{\eta }_{i}^{\eta
_{i}}\right) }{\prod\limits_{i=1}^{n}\overline{\eta }_{i}^{\eta _{i}}}%
=q^{2\eta _{1}+1}t^{-l_{\eta }\left( 1\right) }
\end{equation*}%
has been used. This completes the proof by induction.
\end{proof}

\begin{corollary}
We have 
\begin{equation}
z_{i}E_{\eta }\left( z;q^{-1},t^{-1}\right) =\overline{\eta }%
_{i}\sum\limits _{\substack{ I\in \mathbb{J}_{\eta ,1}^{I},\text{ }i\in I 
\\ c_{I}\left( \eta \right) =\lambda }}\frac{r_{I}^{\left\{ i\right\}
}\left( \overline{c_{I}\left( \eta \right) }\right) k_{\eta }}{%
k_{c_{I}\left( \eta \right) }}E_{c_{I}\left( \eta \right) }\left(
z;q^{-1},t^{-1}\right) .  \label{begin}
\end{equation}
\end{corollary}

\begin{equation*}
\end{equation*}%
We make our final improvement to the formula for $z_{i}E_{\eta }\left(
z;q^{-1},t^{-1}\right) $ by simplifying the coefficient $\overline{\eta }%
_{i}r_{I}^{\left\{ i\right\} }\left( \overline{c_{I}\left( \eta \right) }%
\right) $.

\begin{proposition}
\label{final proposition}Let 
\begin{eqnarray}
\widetilde{B}_{I}\left( z\right)
&:&=\prod\limits_{u=1}^{s}\prod\limits_{j=t_{u-1}+1}^{t_{u}-1}\widehat{b}%
\left( z_{t_{u}},z_{j}\right) \prod\limits_{j=t_{s}+1}^{n}\widehat{b}\left(
qz_{t_{1}},z_{j}\right)  \label{BiTild} \\
&&\times \left( qz_{t_{1}}-t^{-n+1}\right) ,\text{ }t_{0}:=0  \notag
\end{eqnarray}%
and%
\begin{equation*}
\widetilde{\chi }_{I}^{\left\{ i\right\} }\left( z\right) :=\left\{ 
\begin{tabular}{ll}
$\frac{z_{i}}{a\left( z_{i},z_{t_{k+1}}\right) }$ & $;i=t_{k},$ $k=1,...,s-1$
\\ 
$\frac{z_{i}}{a\left( z_{i},qz_{t_{1}}\right) }$ & $;i=t_{s},$%
\end{tabular}%
\right.
\end{equation*}%
where $I=\left\{ t_{1},...,t_{s}\right\} \subseteq \left\{ 1,...,n\right\} ,$
with $1\leq t_{1}<...<t_{s}\leq n$ and $I\neq \emptyset .$ We have 
\begin{equation}
z_{i}E_{\eta }\left( z;q^{-1},t^{-1}\right) =\sum\limits_{\substack{ I\in 
\mathbb{J}_{\eta }^{I},\text{ }i\in I  \\ c_{I}\left( \eta \right) =\lambda 
}}\frac{\widetilde{\chi }_{I}^{\left\{ i\right\} }\left( \overline{\eta }%
\right) A_{I}\left( \overline{\eta }\right) \widetilde{B}_{I}\left( 
\overline{\eta }\right) k_{\eta }}{k_{c_{I}\left( \eta \right) }}%
E_{c_{I}\left( \eta \right) }\left( z;q^{-1},t^{-1}\right) .
\label{final decomposition}
\end{equation}
\end{proposition}

\begin{proof}
It can be seen that for $I$ maximal with respect to $\eta $ we have $%
\widetilde{\chi }_{I}^{\left\{ i\right\} }\left( \overline{\eta }\right) =%
\overline{\eta }_{i}\chi _{I}^{\left\{ i\right\} }\left( \overline{%
c_{I}\left( \eta \right) }\right) $ and $B_{I}\left( \overline{c_{I}\left(
\eta \right) }\right) =\widetilde{B}_{I}\left( \overline{\eta }\right) .$
Since $A_{I}\left( \overline{c_{I}\left( \eta \right) }\right) =A_{I}\left( 
\overline{\eta }\right) ,$ it follows that 
\begin{equation}
\overline{\eta }_{i}r_{I}^{\left\{ i\right\} }\left( \overline{c_{I}\left(
\eta \right) }\right) =\widetilde{\chi }_{I}^{\left\{ i\right\} }\left( 
\overline{\eta }\right) A_{I}\left( \overline{\eta }\right) \widetilde{B}%
_{I}\left( \overline{\eta }\right) .  \label{last substitution}
\end{equation}%
By substituting $\left( \ref{last substitution}\right) $ into $\left( \ref%
{begin}\right) $ we arrive at our final decomposition $\left( \ref{final
decomposition}\right) $.
\end{proof}

\section{The Pieri-type Formula for $r=1$ and the Generalised Binomial
Coefficient$\label{pieri one}$}

The second major result of the paper is to determine the nonsymmetric
analogue of the Pieri-type formula $\left( \ref{pieri}\right) $ for $r=1.$
This formula gives the branching coefficients of Macdonald polynomials of
degree $\left\vert \eta \right\vert +1$ in the expansion of $e_{1}\left(
z\right) =z_{1}+...+z_{n}$ times $E_{\eta }\left( z;q^{-1},t^{-1}\right) $.
These coefficients can be derived as a consequence of Proposition \ref{final
proposition}.

\begin{proposition}
\label{product with elementary function}We have%
\begin{equation}
e_{1}\left( z\right) E_{\eta }\left( z;q^{-1},t^{-1}\right)
=\sum\limits_{I\in \mathbb{J}_{\eta }^{I}}a_{\eta ,c_{I}\left( \eta \right)
}E_{c_{I}\left( \eta \right) }\left( z;q^{-1},t^{-1}\right) {\tiny ,}
\label{elementary product 1}
\end{equation}%
where $a_{\eta ,c_{I}\left( \eta \right) }$ is defined by 
\begin{equation}
a_{\eta ,c_{I}\left( \eta \right) }:=\frac{-\left( q-1\right) d_{\eta
}^{\prime }\left( q^{-1},t^{-1}\right) A_{I}\left( \overline{\eta }\right) 
\widetilde{B}_{I}\left( \overline{\eta }\right) }{q^{\eta _{\min \left(
I\right) }+1}\left( t-1\right) d_{c_{I}(\eta )}^{\prime }\left(
q^{-1},t^{-1}\right) }.  \label{an}
\end{equation}
\end{proposition}

\begin{proof}
Summing (\ref{final decomposition}) over all $i$ and then reversing the
order of summation gives 
\begin{eqnarray}
e_{1}\left( z\right) E_{\eta }\left( z;q^{-1},t^{-1}\right)
&=&\sum\limits_{I\in \mathbb{J}_{\eta }^{I}}\sum\limits_{i\in I}\frac{%
\widetilde{\chi }_{I}^{\left\{ i\right\} }\left( \overline{\eta }\right)
A_{I}\left( \overline{\eta }\right) \widetilde{B}_{I}\left( \overline{\eta }%
\right) k_{\eta }}{k_{c_{I}\left( \eta \right) }}  \label{initial one} \\
&&\times E_{c_{I}\left( \eta \right) }\left( z;q^{-1},t^{-1}\right) .  \notag
\end{eqnarray}%
We have%
\begin{equation}
\sum\limits_{i\in I}\widetilde{\chi }_{I}^{\left\{ i\right\} }\left( 
\overline{\eta }\right) =\frac{\overline{\eta }_{\min \left( I\right)
}\left( 1-q\right) }{\left( t-1\right) }  \label{sum}
\end{equation}%
and%
\begin{equation}
\frac{k_{\eta }}{k_{c_{I}\left( \eta \right) }}=\frac{d_{\eta }^{\prime
}\left( q^{-1},t^{-1}\right) }{q^{2\eta _{\min \left( I\right)
}+1}t^{-l_{\eta }^{\prime }\left( \min \left( I\right) \right)
}d_{c_{I}(\eta )}^{\prime }\left( q^{-1},t^{-1}\right) }.  \label{ratio}
\end{equation}%
Substituting $\left( \ref{sum}\right) $ and $\left( \ref{ratio}\right) $
into $\left( \ref{initial one}\right) $ gives the required result$.$
\end{proof}

\begin{equation*}
\end{equation*}%
On obtaining the Pieri-type formula for $r=1$ we are naturally lead to
deducing an explicit formula for the generalised binomial coefficient $%
\binom{\nu }{\eta }_{q,t}$ when $\left\vert \nu \right\vert =\left\vert \eta
\right\vert +1.$ Generalised binomial coefficients appear in the theory of
Macdonald polynomials. We define nonsymmetric $q-$binomial coefficients $%
\binom{\nu }{\eta }_{q,t}$ according to the generating function formula [\ref%
{peter}]%
\begin{equation}
E_{\eta }\left( z;q^{-1},t^{-1}\right) \prod\limits_{i=1}^{n}\frac{1}{%
\left( z_{i};q\right) _{\infty }}=\sum\limits_{\nu }\binom{\nu }{\eta }%
_{q,t}t^{l\left( \nu \right) -l\left( \eta \right) }\frac{d_{\eta }^{\prime
}\left( q,t\right) }{d_{\nu }^{\prime }\left( q,t\right) }E_{\nu }\left(
z;q^{-1},t^{-1}\right) ,  \label{def binomial}
\end{equation}%
where $\left( z_{i};q\right) _{\infty }$ is the Pockhammer symbol and is
defined as 
\begin{equation}
\left( u;q\right) _{\infty }:=\prod\limits_{j=1}^{\infty }\left(
1-uq^{j-1}\right)  \label{pockhammer}
\end{equation}%
and $l\left( \eta \right) :=\Sigma _{s\in \eta }l_{\eta }\left( s\right) .$
Unlike the classical binomial coefficients 
\begin{equation}
\binom{l}{p}:=\frac{l!}{\left( l-p\right) !p!}  \label{binomial}
\end{equation}%
there is no known explicit formula for $\binom{\nu }{\eta }_{q,t}.$ However,
by restricting our attention to the monomials of degree $1$ in the expansion
of $\prod\nolimits_{i=1}^{n}\frac{1}{\left( z_{i};q\right) _{\infty }}$ we
are able to use Proposition \ref{product with elementary function} to deduce
an explicit formula for $\binom{\nu }{\eta }_{q,t}$ when $\left\vert \nu
\right\vert =\left\vert \eta \right\vert +1$.

\begin{proposition}
\label{binomial formula}Suppose $\left\vert \nu \right\vert =\left\vert \eta
\right\vert +1.$ Then 
\begin{equation}
\binom{\nu }{\eta }_{q,t}=-\frac{A_{I}\left( \overline{\eta }\right) 
\widetilde{B}_{I}\left( \overline{\eta }\right) }{\left( t-1\right) },
\label{binomial coeff}
\end{equation}%
where $\nu =c_{I}\left( \eta \right) .$ If there is no such $I$ such that $%
\nu =c_{I}\left( \eta \right) $ then $\binom{\nu }{\eta }_{q,t}=0.$
\end{proposition}

\begin{proof}
Using $\left( \ref{pockhammer}\right) $ and the identity $\frac{1}{1-u}%
=1+u+u^{2}+...$ we can simplify $\prod\nolimits_{i=1}^{n}\frac{1}{\left(
z_{i};q\right) _{\infty }}$ to 
\begin{equation}
1+\frac{1}{1-q}e_{1}\left( z\right) +\text{ higher order terms}.  \label{1}
\end{equation}%
Equating terms of degree $\left\vert \eta \right\vert +1$ in $\left( \ref%
{def binomial}\right) $ gives 
\begin{equation*}
e_{1}\left( z\right) E_{\eta }\left( z;q^{-1},t^{-1}\right)
=\sum\limits_{\left\vert \nu \right\vert =\left\vert \eta \right\vert +1}%
\binom{\nu }{\eta }_{q,t}\frac{t^{l\left( \nu \right) -l\left( \eta \right)
}\left( 1-q\right) d_{\eta }^{\prime }\left( q,t\right) }{d_{\nu }^{\prime
}\left( q,t\right) }E_{\nu }\left( z;q^{-1},t^{-1}\right) 
\end{equation*}%
(this equation can also be deduced from [\ref{lascoux}, Eq.~(16)]). Comparison with $\left( \ref{an}\right) $ gives 
\begin{equation}
\binom{\nu }{\eta }_{q,t}=\frac{d_{\eta }^{\prime }\left(
q^{-1},t^{-1}\right) }{d_{\eta }^{\prime }\left( q,t\right) }\frac{d_{\nu
}^{\prime }\left( q,t\right) }{d_{\nu }^{\prime }\left( q^{-1},t^{-1}\right) 
}\frac{A_{I}\left( \overline{\eta }\right) \widetilde{B}_{I}\left( \overline{%
\eta }\right) }{q^{\eta _{\min \left( I\right) }+1}\left( t-1\right)
t^{l\left( \nu \right) -l\left( \eta \right) }}.  \label{messy}
\end{equation}%
Since 
\begin{equation*}
\frac{1-x}{1-x^{-1}}=-x,
\end{equation*}%
we have%
\begin{eqnarray}
\frac{d_{\mu }^{\prime }\left( q,t\right) }{d_{\mu }^{\prime }\left(
q^{-1},t^{-1}\right) } &=&\prod\limits_{s\in \mu }(-q^{a_{\mu }\left(
s\right) +1}t^{l_{\mu }\left( s\right) })  \notag \\
&=&\left( -1\right) ^{\left\vert \mu \right\vert }q^{\Sigma _{s\in \mu
}(a_{\mu }\left( s\right) +1)}t^{l\left( \mu \right) }  \label{newby}
\end{eqnarray}%
The final result is obtained by appropriately substituting $\left( \ref%
{newby}\right) $ into $\left( \ref{messy}\right) $ and noting that $\left(
-1\right) ^{\left\vert \nu \right\vert -\left\vert \eta \right\vert }=-1$
while%
\begin{equation*}
\frac{q^{\Sigma _{s\in \nu }(a_{\nu }\left( s\right) +1)}}{q^{\Sigma _{s\in
\eta }(a_{\eta }\left( s\right) +1)}}=q^{\eta _{\min \left( I\right) }+1}.
\end{equation*}
\end{proof}

A viewpoint of the classical binomial coefficients is that they are a ratio
of evaluations of the one variable interpolation polynomial%
\begin{equation*}
f_{p}\left( x\right) :=x\left( x-1\right) ...\left( x-p+1\right) ,
\end{equation*}%
explicitly 
\begin{equation*}
\binom{l}{p}=\frac{f_{p}\left( l\right) }{f_{p}\left( p\right) }.
\end{equation*}%
Similarly in the multivariable nonsymmetric Macdonald polynomial theory the
generalised binomial coefficient $\left( \text{in particular }\left( \ref%
{binomial coeff}\right) \right) $ satisfy [\ref{sahi binomial}] 
\begin{equation}
\binom{\nu }{\eta }_{q,t}=\frac{E_{\eta }^{\ast }\left( \overline{\nu }%
\right) }{E_{\eta }^{\ast }\left( \overline{\eta }\right) }.
\label{binomial evaluation}
\end{equation}%
Comparing $\left( \ref{binomial evaluation}\right) $ with $\left( \ref%
{binomial coeff}\right) $ and making use of $\left( \ref{evaluation}\right) $
gives a new evaluation formula for $E_{\eta }^{\ast }\left( \overline{\nu }%
\right) $ where $\left\vert \nu \right\vert =\left\vert \eta \right\vert +1.$

\begin{corollary}
Suppose $\left\vert \nu \right\vert =\left\vert \eta \right\vert +1.$ Then%
\begin{equation*}
E_{\eta }^{\ast }\left( \overline{\nu }\right) =-\frac{A_{I}\left( \overline{%
\eta }\right) \widetilde{B}_{I}\left( \overline{\eta }\right) }{\left(
t-1\right) }\left( \prod\limits_{i=1}^{n}\overline{\eta }_{i}^{\eta
_{i}}\right) d_{\eta }^{\prime }\left( q^{-1},t^{-1}\right)
\end{equation*}%
where $\nu =c_{I}\left( \eta \right) .$ If there is no such $I$ such that $%
\nu =c_{I}\left( \eta \right) $ then $E_{\eta }^{\ast }\left( \overline{\nu }%
\right) =0$.
\end{corollary}

\section{The Pieri-type Formula for $r=n-1\label{pieri last}$}

In this section we give our last Pieri-type formula, the nonsymmetric
analogue of $\left( \ref{pieri}\right) $ for $r=n-1.$ The result can be
derived almost immediately from the expansion of $e_{1}\left( z\right)
E_{\eta }\left( z\right) $ using the identity [\ref{marshall macdonald}]%
\begin{equation}
E_{\eta }\left( z^{-1};q,t\right) =E_{-\eta ^{R}}\left( z;q,t\right) ,
\label{z inverse}
\end{equation}%
where $\eta ^{R}:=\left( \eta _{n},...,\eta _{1}\right) $.

\begin{proposition}
\label{en-1}Define%
\begin{equation*}
\eta +\left( i^{n}\right) =\left( \eta _{1}+i,...,\eta _{n}+i\right) ,
\end{equation*}%
and%
\begin{equation*}
\eta ^{\prime }:=\eta -\left( \min (\eta )^{n}\right) .
\end{equation*}%
We have 
\begin{equation}
e_{n-1}\left( z\right) E_{\eta }\left( z;q^{-1},t^{-1}\right)
=\sum\limits_{I\in \mathbb{J}_{\nu }^{I}}a_{\nu ,c_{I}\left( \nu \right)
}E_{\lambda +\left( \min (\eta )^{n}\right) }\left( z;q^{-1},t^{-1}\right) ,
\label{last pieri}
\end{equation}%
where $a_{\nu ,c_{I}\left( \nu \right) }$ is defined by $\left( \ref{an}%
\right) ,$ 
\begin{equation}
\nu =\left( -\eta ^{\prime }+\left( \max (\eta ^{\prime })^{n}\right)
\right) ^{R}\text{ and }\lambda =-c_{I}\left( \nu \right) ^{R}+\left( \left(
\max (\nu )+1\right) ^{n}\right) .  \label{definitions}
\end{equation}
\end{proposition}

\begin{proof}
By Proposition \ref{product with elementary function} we have 
\begin{equation*}
e_{1}\left( z\right) E_{\nu }\left( z;q^{-1},t^{-1}\right)
=\sum\limits_{I\in \mathbb{J}_{\nu }^{I}}a_{\nu ,c_{I}\left( \nu \right)
}E_{c_{I}\left( \nu \right) }\left( z;q^{-1},t^{-1}\right) {\tiny .}
\end{equation*}%
Substituting $z$ for $z^{-1}$ and using $\left( \ref{z inverse}\right) $ we
obtain%
\begin{equation*}
e_{1}\left( z^{-1}\right) E_{-\nu ^{R}}\left( z;q^{-1},t^{-1}\right)
=\sum\limits_{I\in \mathbb{J}_{\nu }^{I}}a_{\nu ,c_{I}\left( \nu \right)
}E_{-c_{I}\left( \nu \right) ^{R}}\left( z;q^{-1},t^{-1}\right) .
\end{equation*}%
Multiplying both sides by $z_{1}...z_{n}$ and using the identity $%
z_{1}...z_{n}E_{\eta }\left( z\right) =E_{\eta +\left( 1^{n}\right) }\left(
z\right) $ [\ref{marshall macdonald}] we have%
\begin{equation}
e_{n-1}\left( z\right) E_{-\nu ^{R}}\left( z;q^{-1},t^{-1}\right)
=\sum\limits_{I\in \mathbb{J}_{\nu }^{I}}a_{\nu ,c_{I}\left( \nu \right)
}E_{-c_{I}\left( \nu \right) ^{R}+\left( 1^{n}\right) }\left(
z;q^{-1},t^{-1}\right) .  \label{not yet}
\end{equation}%
Since $\nu =\left( -\eta ^{\prime }+\left( \max (\eta ^{\prime })^{n}\right)
\right) ^{R}$ we have $\eta ^{\prime }=-\nu ^{R}+\left( \max (\nu
)^{n}\right) ,$ and hence, multiplying both sides of $\left( \ref{not yet}%
\right) $ by $\left( z_{1}...z_{n}\right) ^{\max \left( \nu \right) }$ gives 
\begin{equation}
e_{n-1}\left( z\right) E_{\eta ^{\prime }}\left( z;q^{-1},t^{-1}\right)
=\sum\limits_{I\in \mathbb{J}_{\nu }^{I}}a_{\nu ,c_{I}\left( \nu \right)
}E_{\lambda }\left( z;q^{-1},t^{-1}\right) ,  \label{still not yet}
\end{equation}%
where $\lambda $ is defined in $\left( \ref{definitions}\right) $. The final
decomposition $\left( \ref{last pieri}\right) $ is now obtained by
multiplying both sides of $\left( \ref{still not yet}\right) $ by $%
(z_{1}...z_{n})^{\min \left( \eta \right) }.$
\end{proof}

\section{The Classical Limit}

The classical limit in Macdonald polynomial theory refers to setting $%
t=q^{1/\alpha }$ and taking $q\rightarrow 1.$ In particular%
\begin{equation*}
\lim_{t=q^{1/\alpha },\text{ }q\rightarrow 1}E_{\eta }\left( z;q,t\right)
=E_{\eta }\left( z;\alpha \right)
\end{equation*}%
where $E_{\eta }\left( z;\alpha \right) $ is the nonsymmetric Jack
polynomial $\left( \text{for an account of the latter see e.g. [\ref{peter
new}]}\right) $ As remarked in the introduction, the expansion of the
product $e_{1}\left( z\right) E_{\eta }\left( z;\alpha \right) $ in terms of 
$\left\{ E_{\lambda }\left( z;\alpha \right) \right\} $ has been given by
Marshall [\ref{marshall jack}]. We will conclude our study by taking the
classical limit of Proposition \ref{product with elementary function}. First
we recall the result of [\ref{marshall jack}].

\begin{proposition}
\label{marshall prop}We have 
\begin{equation*}
e_{1}\left( z\right) E_{\eta }\left( z;\alpha \right) =\sum\limits_{I\in 
\mathbb{J}_{\eta }^{I}}a_{\eta ,c_{I}\left( \eta \right) }^{\alpha
}E_{c_{I}\left( \eta \right) }\left( z;\alpha \right) ,
\end{equation*}%
where%
\begin{equation}
a_{\eta ,c_{I}\left( \eta \right) }^{\alpha }=\frac{-\alpha ^{2}d_{\alpha
,\eta }^{\prime }A_{\alpha ,I}\left( \frac{\overline{\eta }}{\alpha }\right) 
\widetilde{B}_{\alpha ,I}\left( \frac{\overline{\eta }}{\alpha }\right) }{%
d_{\alpha ,c_{I}\left( \eta \right) }^{\prime }},  \label{new one}
\end{equation}%
The quantities in $\left( \ref{new one}\right) $ are specified by 
\begin{eqnarray}
A_{\alpha ,I}\left( z\right) &:&=a\left( z_{t_{s}}-1,z_{t_{1}}\right)
\prod\limits_{u=1}^{s-1}a\left( z_{t_{u}},z_{t_{u+1}}\right)
\label{Ai alpha} \\
\widetilde{B}_{\alpha ,I}\left( z\right)
&:&=\prod\limits_{u=1}^{s}\prod\limits_{j=t_{u-1}+1}^{t_{u}-1}b\left(
z_{t_{u}},z_{j}\right) \prod\limits_{j=t_{s}+1}^{n}b\left(
z_{t_{1}}+1,z_{j}\right)  \label{BiTild alpha} \\
&&\times \left( z_{t_{1}}+1+\frac{n-1}{\alpha }\right) ,\text{ }t_{0}:=0 
\notag
\end{eqnarray}%
with 
\begin{equation}
a\left( x,y\right) :=\frac{1}{\alpha \left( x-y\right) },\text{ }b\left(
x,y\right) :=\frac{x-y-\frac{1}{\alpha }}{x-y}  \label{a and b alpha}
\end{equation}%
and 
\begin{equation}
d_{\alpha ,\eta }^{\prime }:=\prod\limits_{\left( i,j\right) \in \eta
}\left( \alpha \left( a\left( i,j\right) +1\right) +l\left( i,j\right)
\right) ,  \label{d dash alpha}
\end{equation}%
where $a\left( i,j\right) $ and $l\left( i,j\right) $ are defined by $\left( %
\ref{arm and leg}\right) $ and $I$ by $\left( \ref{I}\right) .$
\end{proposition}

\begin{proof}
Our task is to show that 
\begin{equation}
\lim_{t=q^{1/\alpha },\text{ }q\rightarrow 1}a_{\eta ,c_{I}\left( \eta
\right) }=a_{\eta ,c_{I}\left( \eta \right) }^{\alpha }.
\label{limiting coefficients}
\end{equation}%
Comparing $\left( \ref{d dash}\right) $ with $\left( \ref{d dash alpha}%
\right) ,$ it is immediate that 
\begin{equation*}
\lim_{t=q^{1/\alpha },\text{ }q\rightarrow 1}\frac{\left( q-1\right) d_{\eta
}^{\prime }\left( q^{-1},t^{-1}\right) }{q^{\eta _{\min \left( I\right)
}+1}\left( t-1\right) d_{c_{I}(\eta )}^{\prime }\left( q^{-1},t^{-1}\right) }%
=\alpha ^{2}\frac{d_{\alpha ,\eta }^{\prime }}{d_{\alpha ,c_{I}\left( \eta
\right) }^{\prime }}.
\end{equation*}%
To proceed further, note from $\left( \ref{a and b}\right) $ and $\left( \ref%
{a and b alpha}\right) $ that%
\begin{equation*}
\lim_{t=q^{1/\alpha },\text{ }q\rightarrow 1}\widehat{a}\left(
q^{m}t^{n},q^{m^{\prime }}t^{n^{\prime }}\right) =a\left( \frac{m}{\alpha }%
+n,\frac{m^{\prime }}{\alpha }+n^{\prime }\right)
\end{equation*}%
and%
\begin{equation*}
\lim_{t=q^{1/\alpha },\text{ }q\rightarrow 1}\widehat{b}\left(
q^{m}t^{n},q^{m^{\prime }}t^{n^{\prime }}\right) =b\left( \frac{m}{\alpha }%
+n,\frac{m^{\prime }}{\alpha }+n^{\prime }\right) .
\end{equation*}%
Using this a term-by-term comparison of $\left( \ref{Ai}\right) $ and $%
\left( \ref{BiTild}\right) $ with $\left( \ref{Ai alpha}\right) $ and $%
\left( \ref{BiTild alpha}\right) $ also allows us to conclude that%
\begin{equation*}
\lim_{t=q^{{1}/{\alpha }},\text{ }q\rightarrow 1}A_{I}\left( \overline{%
\eta }\right) =A_{\alpha ,I}\left( \frac{\overline{\eta }}{\alpha }\right)
\end{equation*}%
and 
\begin{eqnarray*}
\lim_{t=q^{{1}/{\alpha }},\text{ }q\rightarrow 1}\widetilde{B}_{I}\left( 
\overline{\eta }\right) &=&\widetilde{B}_{\alpha ,I}\left( \frac{\overline{%
\eta }}{\alpha }\right) . \\
&&
\end{eqnarray*}%
This establishes $\left( \ref{limiting coefficients}\right) ,$ thus
exhibiting Proposition \ref{marshall prop} as a corollary of Proposition \ref%
{product with elementary function}.%
\end{proof}

\subsection*{Acknowledgement.} \textit{I would like to thank my supervisor Peter
Forrester for his guidance and encouragement. This work was supported by an
APA scholarship, the ARC and the University of Melbourne.}


\begin{thebibliography}{99}
\bibitem{peter 2} \label{peter 3}T.H. Baker and P.J. Forrester, A $q-$%
analogue of the Type $A$ Dunkl Operator and Integral Kernal. Int. Math. Res.
Not., 14:667-686, 1997.

\bibitem{cherednik} \label{cherednik} I. Cherednik. Nonsymmetric Macdonald
Polynomials. Int. Math. Res. Not., 10:483-515, 1995.

\bibitem{feigin} \label{japanese}B. Feigin, M. Jimbo, T. Miwa and E. Mukhin,
Symmetric Polynomials Vanishing on the Shifted Diagonals and Macdonald
Polynomials. Int. Math. Res. Not., 18:1015-1034, 2003.

\bibitem{last peter} \label{last peter}P.J. Forrester and E.M. Rains,
Interpretations of Some Parameter Dependent Generalizations of Classical
Matrix Ensembles. Probab. Theory Related Fields, 131:1-61, 2005.

\bibitem{peter} \label{peter}P.J. Forrester, Isomorphisms of Type A Affine
Hecke Algebras and Multivariable Orthogonal Polynomials. Pacific J. Math.,
194:19-41, 1997.

\bibitem{forrester} \label{forrester}P.J. Forrester and D. McAnally,
Pieri-type Formulas for the Nonsymmetric Jack Polynomials. Comment. Math.
Helv. 70:1-24, 2004.

\bibitem{peter new} \label{peter new}P.J. Forrester and T.H. Baker,
Symmetric Jack Polynomials from Nonsymmetric Theory. Annals. Comb.,
3:159-170, 1999.

\bibitem{knop and sahi} \label{knop and sahi}F. Knop and S. Sahi, Difference
Equations and Symmetric Polynomials Defined by their Zeros. Int. Math. Res.
Not., 10:473-486, 1996.

\bibitem{knop} \label{knop}F. Knop, Symmetric and Nonsymmetric Quantum
Capelli Polynomials, Comment. Maths. Helv., 72:84-100, 1997.

\bibitem{lascoux} \label{lascoux}A. Lascoux, Schubert and Macdonald Polynomials, a Parallel. http://phalanstere.univ-mlv.fr/~al/

\bibitem{affine} \label{affine}I. G. Macdonald, Affine Hecke Algebras and
Orthogonal Polynomials. Seminaire Bourbaki, 77,797:1-18, 1995.

\bibitem{macdonald again} \label{macdonald again}I.G. Macdanald, A New Class
of Symmetric Functions. Publ.I.R.M.A, Strasbourg, Actes 20-e Seminaire
Lotharingen, 131-171, 1988.

\bibitem{macdonald} \label{macdonald}I. G. Macdonald, Symmetric Functions
and Hall Polynomials. Oxford University Press, Oxford, 2nd edition, 1995.

\bibitem{marshall macdonald} \label{marshall macdonald}D. Marshall,
Symmetric and Nonsymmetric Macdonald Polynomials. Annals Comb., 3:385-415,
1999.

\bibitem{marshall jack} \label{marshall jack}D. Marshall, The Product of a
Nonsymmetric Jack Polynomial with a Linear Function. Proceedings of the
American Mathematical Society, 131:1817-1827, 2002.

\bibitem{sahi binomial} \label{sahi binomial}S. Sahi, The Binomial Formula
for Nonsymmetric Macdonald Polynomials. q-alg/9703024.

\bibitem{sahi} \label{sahi}S. Sahi, A New Scalar Product for the
Nonsymmetric Jack Polynomials. Int. Math. Res. Not., 20:997-1004, 1996.
\end{thebibliography}
\end{document}